\documentclass[reqno,12pt]{amsart}
\usepackage{amssymb}
\usepackage{amsfonts}
\usepackage{amsmath}
\usepackage{graphicx}
\usepackage{amscd}

\newtheorem{theorem}{Theorem}
\theoremstyle{plain}
\newtheorem{acknowledgement}{Acknowledgement}

\newtheorem{corollary}{Corollary}

\newtheorem{definition}{Definition}

\newtheorem{remark}{Remark}

\numberwithin{equation}{section}

\input{tcilatex}

\begin{document}
\author{}
\title{}
\maketitle

\begin{center}
A \textbf{NEW GENERATING FUNCTION OF (}$q$-\textbf{) BERNSTEIN TYPE
POLYNOMIALS AND THEIR INTERPOLATION FUNCTION}

\bigskip

\bigskip $^{\ast }$Yilmaz\ SIMSEK and $^{\ast \ast }$Mehmet ACIKGOZ
\end{center}

$^{\ast }$University of Akdeniz, Faculty of Arts and Science, Department of
Mathematics,\ 07058-Antalya, TURKEY, ysimsek@akdeniz.edu.tr

$^{\ast \ast }$University of Gaziantep, Faculty of Arts and Science,
Department of Mathematics, 27310 Gaziantep, Turkey, acikgoz@gantep.edu.tr

\begin{center}
{\Large %\textit{\Large{Canada}}
}

\bigskip
\end{center}

\bigskip

\begin{center}
\textbf{{\Large {Abstract}}}
\end{center}

The main object of this paper is to construct a new generating function of
the ($q$-) Bernstein type polynomials. We establish elementary properties of
this function. By using this generating function, we derive recurrence
relation and derivative of the ($q$-) Bernstein type polynomials. We also
give relations between the ($q$-) Bernstein type polynomials, Hermite
polynomials, Bernoulli polynomials of higher-order and\ the second kind
Stirling numbers. By applying Mellin transformation to this generating
function, we define interpolation of the ($q$-) Bernstein type polynomials.
Moreover, we give some applications and questions on approximations of ($q$%
-) Bernstein type polynomials, moments of some distributions in
Statistics.\bigskip

\bigskip


\bigskip

\textbf{CONTENTS}

1. Introduction

2. Preliminary results related to the classical Bernstein, higher-order
Bernoulli and

\ \ \ Hermit polynomials, the second kind Stirling numbers

3. Generating function of the ($q$-) Bernstein type polynomials

4. New identities on ($q$-) Bernstein type polynomials, Hermite polynomials
and

\ \ \ \ first kind Stirling numbers

5. Interpolation function of ($q$-) Bernstein type polynomials

6. Further remarks and observation

\bigskip

\section{Introduction}

In \cite{bernstein}, Bernstein introduced the Bernstein polynomials. Since
that time, many authors have studied on these polynomials and other related
subjects (see cf. \cite{acikgoz}-\cite{wu}), and see also the references
cited in each of these earlier works. The Bernstein polynomials can also be
defined in many different ways. Thus, recently, many applications of these
polynomials have been looked for by many authors. These polynomials have
been used not only for approximations of functions in various areas in
Mathematics but also the other fields such as smoothing in statistics,
numerical analysis and constructing Bezier curve which have many interesting
applications in computer graphics (see cf. \cite{bernstein}, \cite{Guan},
\cite{Oruc}, \cite{SOstrovskaAM}, \cite{Ost}, \cite{Nowak}, \cite{Phillips},
\cite{wu}) and see also the references cited in each of these earlier works.

The ($q$-) Bernstein polynomials have been investigated and studied by many
authors without\textit{\ generating function}. So far, we have not found any
generating function of ($q$-) Bernstein polynomials in the literature.
Therefore, we will consider the following question:

\textit{How can one construct \textbf{generating function} of (}$q$\textit{%
-) Bernstein type polynomials}?

The aim of this paper is to give answer this question and to construct
generating function of the ($q$-) Bernstein type polynomials which is given
in Section 3. By using this generating function, we not only give recurrence
relation and derivative of the ($q$-) Bernstein type polynomials but also
find relations between higher-order Bernoulli polynomials, the second kind
Stirling numbers and the Hermite polynomials. In Section 5, by applying
Mellin transformation to the generating function of the ($q$-) Bernstein
type polynomials, we define interpolation function, which interpolates the ($%
q$-) Bernstein type polynomials at negative integers.

\section{Preliminary results related to the classical Bernstein,
higher-order Bernoulli and Hermit polynomials, the second kind Stirling
numbers}

The Bernstein polynomials play a crucial role in approximation theory and
the other branches of Mathematics and Physics.\ Thus in this section we give
definition and some properties of these polynomials.

Let $f$ be a function on $\left[ 0,1\right] $. The classical Bernstein
polynomials of degree $n$ are defined by%
\begin{equation}
\mathbb{B}_{n}f(x)=\sum_{j=0}^{n}f\left( \frac{j}{n}\right) B_{j,n}(x),\text{
}0\leq x\leq 1,  \label{be1}
\end{equation}%
where $\mathbb{B}_{n}f$\ is called the Bernstein operator and%
\begin{equation}
B_{j,n}(x)=\left(
\begin{array}{c}
n \\
j%
\end{array}%
\right) x^{j}(1-x)^{n-j},  \label{be2}
\end{equation}%
$j=0$, $1$,$\cdots $,$n$ are called the Bernstein basis\ polynomials (or the
Bernstein polynomials of degree $n$). There are $n+1$ $n$th degree Bernstein
polynomials. For mathematical convenience, we set $B_{j,n}(x)=0$ if $j<0$ or
$j>n$ cf. (\cite{bernstein}, \cite{Guan}, \cite{Joy}).

If $f:\left[ 0,1\right] \rightarrow \mathbb{C}$ is a continuous function,
the sequence of Bernstein\ polynomials $\mathbb{B}_{n}f(x)$ converges
uniformly to $f$ on $\left[ 0,1\right] $ cf. \cite{kowalski}.

A recursive definition of the $k$th $n$th Bernstein polynomials can be
written as%
\begin{equation*}
B_{k,n}(x)=(1-x)B_{k,n-1}(x)+xB_{k-1,n-1}(x).
\end{equation*}%
For proof of the above relation see \cite{Joy}.

For $0\leq k\leq n$, derivative of the $n$th degree Bernstein polynomials
are polynomials of degree $n-1$:%
\begin{equation}
\frac{d}{dt}B_{k,n}(t)=n\left( B_{k-1,n-1}(t)-B_{k,n-1}(t)\right) ,
\label{be3}
\end{equation}%
cf. (\cite{bernstein}, \cite{Guan}, \cite{Joy}). On the other hand, in
Section 3, using our a new generating function, we give the other proof of (%
\ref{be3}).

Observe that the Bernstein polynomial of degree $n$, $\mathbb{B}_{n}f$, uses
only the sampled values of $f$ at $t_{nj}=\frac{j}{n}$, $j=0$, $1$,$\cdots $,%
$n$. For $j=0$, $1$,$\cdots $,$n$,%
\begin{equation*}
\beta _{j,n}(x)\equiv (n+1)B_{j,n}(x),\text{ }0\leq x\leq 1,
\end{equation*}%
is the density function of beta distribution $beta(j+1,n+1-j)$.

Let $y_{n}(x)$ be a binomial $b(n,x)$ random variable. Then%
\begin{equation*}
E\left\{ y_{n}(x)\right\} =nt,
\end{equation*}%
and%
\begin{equation*}
var\left\{ y_{n}(x)\right\} =E\left\{ y_{n}(x)-nx\right\} ^{2}=nx(1-x),
\end{equation*}%
and%
\begin{equation*}
\mathbb{B}_{n}f(x)=E\left[ f\left\{ \frac{y_{n}(x)}{n}\right\} \right] ,
\end{equation*}%
cf. \cite{Guan}.

The classical higher-order Bernoulli polynomials $\mathcal{B}_{n}^{(v)}(z)$
defined by means of the following generating function%
\begin{equation}
F^{(v)}(z,t)=e^{tx}\left( \frac{t}{e^{t}-1}\right) ^{v}=\sum_{n=0}^{\infty }%
\mathcal{B}_{n}^{(v)}(z)\frac{t^{n}}{n!}\text{.}  \label{I1}
\end{equation}%
The higher-order Bernoulli polynomials play an important role in the finite
differences and in (analytic) number theory. So, the coefficients in all the
usual cenral-difference formulae for interpolation, numerical
differentiation and integration, and differences in terms of derivatives can
be expressed in terms of these polynomials cf. (\cite{acikgoz}, \cite%
{Norlund}, \cite{LopezTemme}, \cite{SimsekKurt}). These polynomials are
related to the many branches of Mathematics. By substituting $v=1$ into the
above, we have%
\begin{equation*}
F(t)=\frac{te^{tx}}{e^{t}-1}=\sum_{n=1}^{\infty }B_{n}\frac{t^{n}}{n!},
\end{equation*}%
where $B_{n}$ is usual Bernoulli polynomials cf. \cite{simsekJmaa}.

The usual second kind Stirling numbers with parameters $(n,k),$ denote by $%
S(n,k)$, that is the number of partitions of the set $\left\{ 1,2,\cdots
,n\right\} $\ into $k$ non empty set. For any $t$, it is well known that the
second kind Stirling numbers are defined by means of the generating function
cf. (\cite{cakic Milovanovic}, \cite{pinter}, \cite{SimsekARXIV})%
\begin{equation}
F_{S}(t,k)=\frac{(-1)^{k}}{k!}(1-e^{t})^{k}=\sum_{n=0}^{\infty }S(n,k)\frac{%
t^{n}}{n!}.  \label{I2}
\end{equation}%
These numbers play an important role in many branches of Mathematics, for
example, combinatorics, number theory, discrete probability distributions
for finding higher order moments. In \cite{Joarder}, Joarder and \ Mahmood
demonstrated the application of Stirling numbers of the second kind in
calculating moments of some discrete distributions, which are binomial
distribution, geometric distribution and negative binomial distribution.

The Hermite polynomials defined by the following generating function:

For $z$, $t\in \mathbb{C}$,%
\begin{equation}
e^{2zt-t^{2}}=\sum_{n=0}^{\infty }H_{n}(z)\frac{t^{n}}{n!},  \label{I3}
\end{equation}%
which gives the Cauchy-type integral%
\begin{equation*}
H_{n}(z)=\frac{n!}{2\pi i}\int_{\mathcal{C}}e^{2zt-t^{2}}\frac{dt}{t^{n+1}},
\end{equation*}%
where $\mathcal{C}$ is a circle around the origin and the integration is in
positive direction cf. \cite{LopezTemme}. The Hermite polynomials play a
crucial role in certain limits of the classical orthogonal polynomials.
These polynomials are related to the higher-order Bernoulli polynomials,
Gegenbauer polynomials, Laguerre polynomials, the Tricomi-Carlitz
polynomials and Buchholz polynomials, cf. \cite{LopezTemme}. These
polynomials also play a crucial role in not only in Mathematics but also in
Physics and in the other sciences. In section 4 we give relation between the
Hermite polynomials and ($q$-) Bernstein type polynomials.

\section{Generating Function of the Bernstein type polynomials}

Let $\left\{ B_{k,n}(x)\right\} _{0\leq k\leq n}$ be a sequence of Bernstein
polynomials. The aim of this section is to construct generating function of
the sequence $\left\{ B_{k,n}(x)\right\} _{0\leq k\leq n}$. It is well known
that most of generating functions are obtained from the recurrence formulae.
However, we do not use the recurrence formula of the Bernstein polynomials
for constructing generating function of them.

We now give the following notation:%
\begin{equation*}
\lbrack x]=[x:q]=\left\{
\begin{array}{c}
\frac{1-q^{x}}{1-q}\text{, }q\neq 1 \\
\\
x\text{, }q=1.%
\end{array}%
\right.
\end{equation*}

If $q\in \mathbb{C}$, we assume that $\mid q\mid <1$.

We define%
\begin{eqnarray}
F_{k,q}(t,x) &=&(-1)^{k}t^{k}\exp \left( \left[ 1-x\right] t\right)
\label{F1} \\
&&\times \sum_{m,l=0}\left(
\begin{array}{c}
k+l-1 \\
l%
\end{array}%
\right) \frac{q^{l}S(m,k)\left( x\log q\right) ^{m}}{m!},  \notag
\end{eqnarray}%
where $\left\vert q\right\vert <1$, $\exp (x)=e^{x}$\ and $S(m,k)$ denotes
the second kind Stirling numbers and%
\begin{equation*}
\sum_{m,l=0}f(m)g(l)=\sum_{m=0}^{\infty }f(m)\sum_{l=0}^{\infty }g(l).
\end{equation*}%
By (\ref{F1}), we define the following a new generating function of
polynomial $Y_{n}(k;x;q)$ by%
\begin{equation}
F_{k,q}(t,x)=\sum_{n=k}^{\infty }Y_{n}(k;x;q)\frac{t^{n}}{n!},  \label{F2}
\end{equation}%
where $t\in \mathbb{C}$.

Observe that if $q\rightarrow 1$ in (\ref{F2}), we have%
\begin{equation*}
Y_{n}(k;x;q)\rightarrow B_{k,n}(x),
\end{equation*}%
hence%
\begin{equation*}
F_{k}(t,x)=\sum_{n=k}^{\infty }B_{k,n}(x)\frac{t^{n}}{n!}.
\end{equation*}

From (\ref{F2}), we obtain the following theorem.

\begin{theorem}
Let $n$ be a positive integer with $k\leq n$. Then we have%
\begin{eqnarray}
Y_{n}(k;x;q) &=&\left(
\begin{array}{c}
n \\
k%
\end{array}%
\right) \frac{(-1)^{k}k!}{(1-q)^{n-k}}  \label{F2a} \\
&&\times \sum_{m,l=0}\sum_{j=0}^{n-k}\left(
\begin{array}{c}
k+l-1 \\
l%
\end{array}%
\right) \left(
\begin{array}{c}
n-k \\
k%
\end{array}%
\right) \frac{(-1)^{j}q^{l+j(1-x)}S(m,k)\left( x\log q\right) ^{m}}{m!}.
\notag
\end{eqnarray}
\end{theorem}

By using (\ref{F1}) and (\ref{F2}), we obtain%
\begin{eqnarray}
F_{k,q}(t,x) &=&\frac{\left( \left[ x\right] t\right) ^{k}}{k!}\exp (\left[
1-x\right] t)  \label{be33} \\
&=&\sum_{n=k}^{\infty }Y_{n}(k;x;q)\frac{t^{n}}{n!}.  \notag
\end{eqnarray}%
The generating function $F_{k,q}(t,x)$ depends on integer parameter $k$,
real variable $x$ and complex variable $q$ and $t$. Therefore the
proprieties of this function are closely related to these variables and
parameter. By using this function, we give many properties of the ($q$-)
Bernstein type polynomials and the other well-known special numbers and
polynomials. By applying Mellin transformation to this function, in Section
5, we construct interpolation function of the ($q$-) Bernstein type
polynomials.

By the \textit{umbral calculus }convention in (\ref{be33}), then we obtain%
\begin{equation}
\frac{\left( \left[ x\right] t\right) ^{k}}{k!}\exp (\left[ 1-x\right]
t)=\exp \left( Y(k;x;q)t\right) .  \label{be3Yn}
\end{equation}%
By using the above, we obtain all recurrence formulae of $Y_{n}(k;x;q)$\ as
follows:%
\begin{equation*}
\frac{\left( \left[ x\right] t\right) ^{k}}{k!}=\sum_{n=0}^{\infty }\left(
Y(k;x;q)-\left[ 1-x\right] \right) ^{n}\frac{t^{n}}{n!},
\end{equation*}%
where each occurrence of $Y^{n}(k;x;q)$ by $Y_{n}(k;x;q)$ (symbolically $%
Y^{n}(k;x;q)\rightarrow Y_{n}(k;x;q)$).

By (\ref{be3Yn}),%
\begin{equation*}
\left[ u+v\right] =\left[ u\right] +q^{u}\left[ v\right]
\end{equation*}%
and%
\begin{equation*}
\left[ -u\right] =-q^{u}\left[ u\right] ,
\end{equation*}%
we obtain the following corollary:

\begin{corollary}
\label{corl-1} Let $n$ be a positive integer with $k\leq n$. Then we have%
\begin{equation*}
Y_{n+k}(k;x;q)=\left(
\begin{array}{c}
n+k \\
k%
\end{array}%
\right) \sum_{j=0}^{n}(-1)^{j}q^{j(1-x)}\left[ x\right] ^{j+k}.
\end{equation*}
\end{corollary}

\begin{remark}
By Corollary \ref{corl-1}, for all $k$ with $0\leq k\leq n$, we see that%
\begin{equation*}
Y_{n+k}(k;x;q)=\left(
\begin{array}{c}
n+k \\
k%
\end{array}%
\right) \sum_{j=0}^{n}(-1)^{j}q^{j(1-x)}\left[ x\right] ^{j+k},
\end{equation*}%
or%
\begin{equation*}
Y_{n+k}(k;x;q)=\left(
\begin{array}{c}
n+k \\
k%
\end{array}%
\right) \left[ x\right] ^{k}\left[ 1-x\right] ^{n}.
\end{equation*}%
The polynomials $Y_{n+k}(k;x;q)$ are so-called $q$-\textbf{Bernstein-type
polynomials}. It is easily seen that%
\begin{equation*}
\lim_{q\rightarrow 1}Y_{n+k}(k;x;q)=B_{k,n+k}(x)=\left(
\begin{array}{c}
n+k \\
k%
\end{array}%
\right) x^{k}\left( 1-x\right) ^{n},
\end{equation*}%
which give us (\ref{be2}).
\end{remark}

By using derivative operator%
\begin{equation*}
\frac{d}{dx}\left( \lim_{q\rightarrow 1}Y_{n+k}(k;x;q)\right)
\end{equation*}%
in (\ref{F1}), we obtain%
\begin{eqnarray*}
&&\sum_{n=k}^{\infty }\frac{d}{dx}\left( Y_{n}(k;x;1)\right) \frac{t^{n}}{n!}
\\
&=&\sum_{n=k}^{\infty }nY_{n-1}(k-1;x;1)\frac{t^{n}}{n!}-\sum_{n=k}^{\infty
}nY_{n-1}(k;x;1)\frac{t^{n}}{n!}.
\end{eqnarray*}%
Consequently, we have%
\begin{equation*}
\frac{d}{dx}\left( Y_{n}(k;x;1)\right) =nY_{n-1}(k-1;x;1)-nY_{n-1}(k;x;1),
\end{equation*}%
or%
\begin{equation*}
\frac{d}{dx}\left( B_{k,n}(x)\right) =nB_{k-1,n-1}(x)-nB_{k,n-1}(x).
\end{equation*}

Observe that by using our generating function, we give different proof of (%
\ref{be3}).

Let $f$ be a function on $\left[ 0,1\right] $. The ($q$-) Bernstein type
polynomial of degree $n$ is defined by%
\begin{equation*}
\mathbb{Y}_{n}f(x)=\sum_{j=0}^{n}f\left( \frac{\left[ j\right] }{\left[ n%
\right] }\right) Y_{n}(j;x;q),
\end{equation*}%
where $0\leq x\leq 1$. $\mathbb{Y}_{n}$ is called the ($q$-) Bernstein type
operator and $Y_{n}(j;x;q)$, $j=0,\cdots ,n$, defined in (\ref{F2a}), are
called the ($q$-) Bernstein type (basis) polynomials.

\section{New identities on Bernstein type polynomials, Hermite polynomials
and\ first kind Stirling numbers}

\begin{theorem}
Let $n$ be a positive integer with $k\leq n$. Then we have%
\begin{equation*}
Y_{n}(k;x;q)=\left[ x\right] ^{k}\sum_{j=0}^{n}\left(
\begin{array}{c}
n \\
j%
\end{array}%
\right) \mathcal{B}_{j}^{(k)}\left( \left[ 1-x\right] \right) S(n-j,k),
\end{equation*}%
where $\mathcal{B}_{j}^{(k)}(x)$ and $S(n,k)$ denote the classical
higher-order Bernoulli polynomials and the second kind Stirling numbers,
respectively.
\end{theorem}

\begin{proof}
By using (\ref{I1}), (\ref{I2}) and (\ref{F2}), we obtain%
\begin{eqnarray*}
&&\sum_{n=k}^{\infty }Y_{n}(k;x;q)\frac{t^{n}}{n!} \\
&=&\left[ x\right] ^{k}\sum_{n=0}^{\infty }S(n,k)\frac{t^{n}}{n!}%
\sum_{n=0}^{\infty }\mathcal{B}_{j}^{(k)}\left( \left[ 1-x\right] \right)
\frac{t^{n}}{n!}.
\end{eqnarray*}%
By using Cauchy product in the above, we have%
\begin{eqnarray*}
&&\sum_{n=k}^{\infty }Y(k,n;x;q)\frac{t^{n}}{n!} \\
&=&\left[ x\right] ^{k}\sum_{n=0}^{\infty }\sum_{j=0}^{n}\mathcal{B}%
_{j}^{(k)}\left( \left[ 1-x\right] \right) S(n-j,k)\frac{t^{n}}{j!\left(
n-j\right) !}.
\end{eqnarray*}%
From the above, we have%
\begin{eqnarray}
&&\sum_{n=k}^{\infty }Y_{n}(k;x;q)\frac{t^{n}}{n!}  \label{b1} \\
&=&\left[ x\right] ^{k}\sum_{n=0}^{k-1}\sum_{j=0}^{n}\mathcal{B}%
_{j}^{(k)}\left( \left[ 1-x\right] \right) S(n-j,k)\frac{t^{n}}{j!\left(
n-j\right) !}  \notag \\
&&+\left[ x\right] ^{k}\sum_{n=k}^{\infty }\sum_{j=0}^{n}\mathcal{B}%
_{j}^{(k)}\left( \left[ 1-x\right] \right) S(n-j,k)\frac{t^{n}}{j!\left(
n-j\right) !}.  \notag
\end{eqnarray}%
By comparing coefficients of $t^{n}$ in the both sides of the above
equation, we arrive at the desired result.
\end{proof}

\begin{remark}
In \cite{H. W. Gould}, Gould gave a different relation between the Bernstein
polynomials, generalized Bernoulli polynomials and the second kind Stirling
numbers. Oruc and Tuncer \cite{Oruc} gave relation between the $q$-Bernstein
polynomials and the second kind $q$-Stirling numbers. In \cite{Nowak}, Nowak
studied on approximation properties for generalized $q$-Bernstein
polynomials and also obtained Stancu operators or Phillips polynomials.
\end{remark}

From (\ref{b1}), we get the following corollary:

\begin{corollary}
Let $n$ be a positive integer with $k\leq n$. Then we have%
\begin{equation*}
\left[ x\right] ^{k}\sum_{n=0}^{k-1}\sum_{j=0}^{n}\frac{\mathcal{B}%
_{j}^{(k)}\left( \left[ 1-x\right] \right) S(n-j,k)}{j!\left( n-j\right) !}%
=0.
\end{equation*}
\end{corollary}

\begin{theorem}
Let $n$ be a positive integer with $k\leq n$. Then we have%
\begin{equation*}
H_{n}(1-y)=\frac{k!}{y^{k}}\sum_{n=0}^{\infty }Y_{n+k}(k;y;q)\frac{2^{n}}{%
\left( n+k\right) !}.
\end{equation*}
\end{theorem}

\begin{proof}
By (\ref{I3}), we have%
\begin{equation*}
e^{2zt}=\sum_{n=0}^{\infty }\frac{t^{2n}}{n!}\sum_{n=0}^{\infty }H_{n}(z)%
\frac{t^{n}}{n!}.
\end{equation*}%
By Cauchy product in the above, we obtain%
\begin{equation}
e^{2zt}=\sum_{n=0}^{\infty }\left( \sum_{j=0}^{n}\left(
\begin{array}{c}
n \\
j%
\end{array}%
\right) H_{j}(z)\right) \frac{t^{2n-j}}{n!}.  \label{b2}
\end{equation}%
By substituting $z=1-y$ into (\ref{b2}), we have%
\begin{equation*}
\sum_{n=0}^{\infty }\left( \sum_{j=0}^{n}\left(
\begin{array}{c}
n \\
j%
\end{array}%
\right) H_{j}(1-y)\right) \frac{t^{2n-j}}{n!}=\frac{k!}{y^{k}}%
\sum_{n=0}^{\infty }\left( 2^{n}Y_{n+k}(k;y;q)\right) \frac{t^{n}}{\left(
n+k\right) !}.
\end{equation*}%
By comparing coefficients of $t^{n}$ in the both sides of the above
equation, we arrive at the desired result.
\end{proof}

\section{Interpolation Function of the ($q$-) Bernstein type polynomials}

The classical Bernoulli numbers interpolate by Riemann zeta function, which
has a profound effect on number theory and complex analysis. Thus, we
construct interpolation function of the ($q$-) Bernstein type polynomials.

For $z\in \mathbb{C}$, and $x\neq 1$, by applying the Mellin transformation
to (\ref{F1}), we get%
\begin{equation*}
S_{q}(z,k;x)=\frac{1}{\Gamma (s)}\int_{0}^{\infty }t^{z-k-1}F_{k,q}(-t,x)dt.
\end{equation*}%
By using the above equation, we defined interpolation function of the
polynomials, $Y_{n}(k;x;q)$ as follows:

\begin{definition}
Let $z\in \mathbb{C}$ and $x\neq 1$. We define%
\begin{equation}
S_{q}(z,k;x)=\left( 1-q\right) ^{z-k}\sum_{m,l=0}\left(
\begin{array}{c}
z+l-1 \\
l%
\end{array}%
\right) \frac{q^{l(1-x)}S(m,k)\left( x\log q\right) ^{m}}{m!}.  \label{F1A}
\end{equation}
\end{definition}

By using (\ref{F1A}), we obtain%
\begin{equation*}
S_{q}(z,k;x)=\frac{(-1)^{k}}{k!}\left[ x\right] ^{k}\left[ 1-x\right] ^{-z},
\end{equation*}%
where $z\in \mathbb{C}$ and $x\neq 1$.

By (\ref{F1A}), we have $S_{q}(z,k;x)\rightarrow S(z,k;x)$ as $q\rightarrow
1 $. Thus we have%
\begin{equation*}
S(z,k;x)=\frac{(-1)^{k}}{k!}x^{k}\left( 1-x\right) ^{-z}.
\end{equation*}%
By substituting $x=1$ into the above, we have%
\begin{equation*}
S(z,k;1)=\infty .
\end{equation*}%
We now evaluate the $m$th $z$-derivatives of $S(z,k;x)$\ as follows:%
\begin{equation}
\frac{\partial ^{m}}{\partial z^{m}}S(z,k;x)=\log ^{m}\left( \frac{1}{1-x}%
\right) S(z,k;x),  \label{F1b}
\end{equation}%
where $x\neq 1$.

By substituting $z=-n$ into (\ref{F1A}), we obtain%
\begin{equation*}
S_{q}(-n,k;x)=\frac{1}{\left( 1-q\right) ^{n}}\sum_{m,l=0}\left(
\begin{array}{c}
-n+l-1 \\
l%
\end{array}%
\right) \frac{q^{l(1-x)}S(m,k)\left( x\log q\right) ^{m}}{m!}.
\end{equation*}%
By substituting (\ref{F2a}) into the above, we arrive at the following
theorem, which relates the polynomials $Y_{n+k}(k;x;q)$ and the function $%
S_{q}(z,k;x)$.

\begin{theorem}
Let $n$ be a positive integer with $k\leq n$ and $0<x<1$. Then we have%
\begin{equation*}
S_{q}(-n,k;x)=\frac{(-1)^{k}n!}{(n+k)!}Y_{n+k}(k;x;q).
\end{equation*}
\end{theorem}

\begin{remark}
\begin{eqnarray*}
\lim_{q\rightarrow 1}S_{q}(-n,k;x) &=&S(-n,k;x) \\
&=&\frac{(-1)^{k}n!}{(n+k)!}x^{k}\left( 1-x\right) ^{n} \\
&=&\frac{(-1)^{k}n!}{(n+k)!}B_{k,n+k}(x).
\end{eqnarray*}%
Therefore, for $0<x<1$, the function%
\begin{equation*}
S(z,k;x)=\frac{(-1)^{k}}{k!}x^{k}\left( 1-x\right) ^{-z}
\end{equation*}%
interpolates the classical Bernstein polynomials of degree $n$ at negative
integers.
\end{remark}

By substituting $z=-n$ into (\ref{F1b}), we obtain the following corollary.

\begin{corollary}
Let $n$ be a positive integer with $k\leq n$ and $0<x<1$. Then we have%
\begin{equation*}
\frac{\partial ^{m}}{\partial z^{m}}S(-n,k;x)=\frac{(-1)^{k}n!}{(n+k)!}%
B_{k,n+k}(x)\log ^{m}\left( \frac{1}{1-x}\right) .
\end{equation*}
\end{corollary}

\section{Further remarks and observation}

The Bernstein polynomials are used for important applications in many
branches of Mathematics and the other sciences, for instance, approximation
theory, probability theory, statistic theory, number theory, the solution of
the differential equations, numerical analysis, constructing Bezier curve, $%
q $-analysis, operator theory and applications in computer graphics. Thus we
look for the applications of our new functions and the ($q$-) Bernstein type
polynomials. Due to Oruc and Tuncer \cite{Oruc}, the $q$-Bernstein
polynomials shares the well-known shape-preserving properties of the
classical Bernstein polynomials. When the function $f$ is convex then%
\begin{equation*}
\beta _{n-1}(f,x)\geq \beta _{n}(f,x)\text{ for }n>1\text{ and }0<q\leq 1,
\end{equation*}%
where%
\begin{equation*}
\beta _{n}(f,x)=\sum_{r=0}^{n}f_{r}\left[
\begin{array}{c}
n \\
r%
\end{array}%
\right] x^{r}\prod_{s=0}^{n-r-1}\left( 1-q^{s}x\right)
\end{equation*}%
and%
\begin{equation*}
\left[
\begin{array}{c}
n \\
r%
\end{array}%
\right] =\frac{\left[ n\right] \cdots \left[ n-r+1\right] }{\left[ r\right] !%
}.
\end{equation*}%
As a consequence of this one can show that the approximation to convex
function by the $q$-Bernstein polynomials is one sided, with $\beta
_{n}f\geq f$ for all $n$. $\beta _{n}f$ behaves is very nice way when one
vary the parameter $q$. In \cite{acikgoz}, the authors gave some
applications on the approximation theory related to Bernoulli and Euler
polynomials.

We conclude this section by the following questions:

1) \textit{How can one demonstrate approximation by (}$q$-) \textit{%
Bernstein type polynomials, }$Y_{n+k}(k;x;q)$?

2) \textit{Is it possible to define uniform expansions of the (}$q$-)
\textit{Bernstein type polynomials, }$Y_{n+k}(k;x;q)$\textit{?}

3) \textit{Is it possible to give applications of the (}$q$\textit{-)
Bernstein type polynomials in calculating moments of some distributions in
Statistics, }$Y_{n+k}(k;x;q)$\textit{?}

4) \textit{How can one give relations between the (}$q$\textit{-) Bernstein
type polynomials, }$Y_{n+k}(k;x;q)$\textit{\ and the Milnor algebras.}

\begin{acknowledgement}
The first author is supported by the research fund of Akdeniz University.
\end{acknowledgement}

\end{document}